\documentclass[11pt]{amsart}
\def\a{\alpha}

\usepackage{pst-all}
\usepackage{amsmath}
\usepackage[latin1]{inputenc}
\usepackage{amsfonts}
\usepackage{amssymb}
\usepackage{multicol}
\usepackage{verbatim}
\usepackage{amsthm}
\usepackage{amscd}
\usepackage{pb-diagram}
\usepackage[mathscr]{eucal}
\usepackage[latin1]{inputenc}

\usepackage{amsmath}

\usepackage{amsfonts}

\usepackage{amssymb}

\usepackage{amsthm}

\usepackage{graphics}

\newcommand{\ZZ}{\mathbb{Z}}

\newcommand{\CC}{\mathbb{C}}

\newcommand{\NN}{\mathbb{N}}

\newcommand{\QQ}{\mathbb{Q}}

\newcommand{\Glie}{\mathfrak{g}}

\newcommand{\Yim}{\mathcal{Y}}

\newcommand{\Hlie}{\mathfrak{h}}

\newcommand{\demo}{\noindent {\it \small Proof:}\quad}

\newcommand{\U}{\mathcal{U}}

\newtheorem{thm}{Theorem}[section]

\newtheorem{defi}[thm]{Definition}

\newtheorem{cor}[thm]{Corollary}

\newtheorem{prop}[thm]{Proposition}

\newtheorem{lem}[thm]{Lemma}

\newtheorem{rem}[thm]{Remark}

\newtheorem{ex}[thm]{Example}

\usepackage[all]{xy}
\title{Simple tensor products}

\author[David Hernandez]{David Hernandez$^1$}\thanks{$^1$Supported partially by ANR through Project "G\'eom\'etrie et Structures
  Alg\'ebriques Quantiques"}
\address{CNRS, \'Ecole Polytechnique and \'Ecole Normale Sup\'erieure, 45, rue d'Ulm 75005 Paris, FRANCE}
\email{dhernand@dma.ens.fr}

\begin{document}

\begin{abstract} Let $\mathcal{F}$ be the category of finite-dimensional representations of an arbitrary quantum affine algebra. We prove that a tensor product $S_1\otimes \cdots \otimes S_N$ of simple objects of $\mathcal{F}$ is simple if and only $S_i\otimes S_j$ is simple for any $i < j$.
\end{abstract}

\maketitle

\tableofcontents

\section{Introduction}

Let $q\in\CC^*$ which is not a root of unity and let $\U_q(\Glie)$ be a quantum affine algebra (not necessarily simply-laced or untwisted). Let $\mathcal{F}$ be the tensor category of finite-dimensional representations of $\U_q(\Glie)$. We prove the following result, expected in various papers of the vast literature about $\mathcal{F}$.

\begin{thm}\label{factg} Let $S_1,\cdots, S_N$ be simple objects of $\mathcal{F}$. The tensor product $S_1\otimes \cdots \otimes S_N$ is simple if and only if $S_i\otimes S_j$ is simple for any $i < j$. 
\end{thm}

The "only if" part of the statement is known : it is an immediate consequence of the commutativity of the Grothendieck ring $\text{Rep}(\U_q(\Glie))$ of $\mathcal{F}$ proved in \cite{Fre} (see \cite{h8} for the twisted types). This will be explained in more details in Section \ref{quatre}. Note that the condition $i < j$ can be replaced by $i\neq j$. Indeed, although in general the two modules $S_i\otimes S_j$ and $S_j\otimes S_i$ are not isomorphic, they are isomorphic if one of them is simple. The "if" part of the statement is the main result of this paper.

If the reader is not familiar with the representation theory of quantum affine algebras, he may wonder why such a result is non trivial. Indeed, in tensor categories associated to "classical" representation theory, there are "few" non trivial tensor products of representations which are simple. For instance, let $V, V'$ be non-zero simple finite-dimensional modules of a simple algebraic group in characteristic $0$. Then, it is well-known that $V\otimes V'$ is simple if and only if $V$ or $V'$ is of dimension $1$. But in positive characteristic there are examples of non trivial simple tensor products given by the Steinberg theorem \cite{s}. And in $\mathcal{F}$ there are "many" simple tensor products of non trivial simple representations. For instance, it is proved in \cite{cp} that for $\Glie = \hat{sl}_2$ an arbitrary simple object $V$ of $\mathcal{F}$ is {\it real}, i.e. $V\otimes V$ is simple. Although it is known \cite{le} that there are non real simple objects in $\mathcal{F}$ when $\Glie$ is arbitrary, many other examples of non trivial simple tensor products can be found in \cite{hl}.

The statement of Theorem \ref{factg} has been conjectured and proved by several authors in various special cases.

\begin{itemize}

\item The result is proved for $\Glie = \hat{sl}_2$ in \cite{cp, Cha2}. 

\item A similar result is proved for a special class of modules of the Yangian of $gl_n$ attached to skew Young diagrams in \cite{NT}. 

\item The result is proved for tensor products of fundamental representations in \cite{ak, Fre2}. 

\item The result is proved for a special class of tensor products satisfying an irreducibility criterion in \cite{Chari2} for the untwisted types.

\item The result is proved for a ``small'' subtensor category of $\mathcal{F}$ when $\Glie$ is simply-laced in \cite{hl}. 

\end{itemize}

\noindent So, even in the case $\Glie = \hat{sl}_3$, the result had not been established. Our complete proof is valid for arbitrary simple objects of $\mathcal{F}$ and for arbitrary $\Glie$.

\noindent Note that the statement of Theorem \ref{factg} allows to produce simple tensor products $V\otimes V'$ where $V = S_1\otimes \cdots \otimes S_k$ and $V' = S_{k+1}\otimes \cdots \otimes S_N$. Besides it implies that $S_1\otimes \cdots \otimes S_N$ is real if we assume that the $S_i$ are real in addition to the assumptions of Theorem \ref{factg}.

Our result is stated in terms of the tensor structure of $\mathcal{F}$. Thus, it is purely representation theoretical. But we have three additional motivations, related respectively to physics, topology, combinatorics, and also to other structures of $\mathcal{F}$.

First, although the category $\mathcal{F}$ is not braided (in general $V\otimes V'$ is not isomorphic to $V'\otimes V$), $\U_q(\Glie)$ has a {\it universal $R$-matrix} in a completion of the tensor product $\U_q(\Glie)\otimes \U_q(\Glie)$. In general the universal $R$-matrix can not be specialized to finite-dimensional representations, but it gives rise to {\it intertwining operators} $V(z)\otimes V'\rightarrow V' \otimes V(z)$ which depend meromorphically on a formal parameter $z$ (see \cite{Fre0, ks}; here the representation $V(z)$ is obtained by homothety of spectral parameter). From the physical point of view, it is an important question to localize the zeros and poles of these operators.
The reducibility of tensor products of objects in $\mathcal{F}$ is known to have strong relations with this question. This is the first motivation to study irreducibility of tensor products in terms of irreducibility of tensor products of pairs of constituents (see \cite{ak} for instance). 

Secondly, if $V\otimes V'$ is simple the universal $R$-matrix can be specialized and we get a well-defined intertwining operator $V\otimes V'\rightarrow V'\otimes V$. In general the action of the $R$-matrix is not trivial (see examples in \cite{jm2}). As the $R$-matrix satisfies the Yang-Baxter equation, when $V$ is real we can define an action of the braid group $\mathcal{B}_N$ on $V^{\otimes N}$ (as for representations of quantum groups of finite type). It is known \cite{rt} that such situations are important to construct {\it topological invariants}.

Finally, in a tensor category, there are natural important questions such as the parametrization of simple objects or the decomposition of tensor products of simple objects in the Grothendieck ring. But another problem of the same importance is the factorization of simple objects $V$ in {\it prime objects}, i.e. the decomposition $V = V_1\otimes \cdots \otimes V_N$ where the $V_i$ can not be written as a tensor product of non trivial simple objects.
This problem for $\mathcal{F}$ is one of the main motivation in \cite{hl}.
When we have established that the tensor products of some pairs of prime representations are simple, Theorem \ref{factg} gives the factorization of arbitrary tensor products of these representations. 
This factorization problem is related to the program of realization of {\it cluster algebras} in $\text{Rep}(\U_q(\Glie))$ initiated in \cite{hl} when $\Glie$ is simply-laced (see more results in this direction in \cite{Nacl}). Cluster algebras have a distinguished set of generators called {\it cluster variables}, and (in finite cluster type) a distinguished linear basis of certain products of cluster variables called {\it cluster monomials}. In the general framework of {\it monoidal categorification} of cluster algebras \cite{hl}, the cluster monomials should correspond to simple modules. 
Theorem \ref{factg} reduces the proof of the irreducibility of tensor products of representations corresponding to cluster variables to the proof of the irreducibility of the tensor products of pairs of simple representations corresponding to cluster variables. 
To conclude with the motivations, Theorem \ref{factg} will be used in the future to establish monoidal categorifications associated to non necessarily simply-laced quantum affine algebras, involving categories different than the small subcategories considered in \cite{hl, Nacl}.

The paper is organized as follows. In Section \ref{un} we give reminders on the category $\mathcal{F}$, in particular on $q$-characters (which will be one of the main tools for the proof). In Section \ref{deux} we prove a general result about tensor products of $l$-weight vectors. In Section \ref{bis} we reduce the problem.
In Section \ref{trois} we introduce upper, lower $q$-characters and we prove several formulae for them. In Section \ref{quatre} we end the proof of the main Theorem \ref{factg}. In Section \ref{cinq} we give some final comments.

{\bf Acknowledgments :} The author is very grateful to Bernard Leclerc for having encouraged him to prove this conjecture. He would like to thank Michio Jimbo and Jean-Pierre Serre for their comments and the Newton Institute in Cambridge where this work was finalized.

\section{Finite-dimensional representations of quantum affine algebras}\label{un}

We recall the main definitions and the main properties of the finite-dimensional representations of quantum affine algebras. For more details, we refer to \cite{Cha2, ch} (untwisted types) and to \cite{Cha5, h8} (twisted types).

\subsection{} In this subsection we shall give all definitions which are sufficient to state Theorem \ref{factg}. All vector spaces, algebras and tensor products are defined over $\CC$.

Fix $h\in\CC$ satisfying $q = e^h$. Then $q^r = e^{hr}$ is well-defined for any $r\in\QQ$.

Let $C = (C_{i,j})_{0\leq i,j\leq n}$ be a {\it generalized Cartan matrix} \cite{kac}, i.e. for $0\leq i,j\leq n$
we have $C_{i,j}\in\ZZ$, $C_{i,i} = 2$, and for $0\leq i\neq j\leq n$ we have $C_{i,j} \leq 0$, $(C_{i,j} = 0 \Leftrightarrow C_{j,i} = 0)$. We suppose that $C$ is {\it indecomposable}, i.e. there is no proper $J\subset \{0,\cdots, n\}$ such that $C_{i,j} = 0$ for any $(i,j)\in J\times (\{0,\cdots, n\}\setminus J)$. Moreover we suppose that $C$ is of {\it affine type}, i.e. all proper principal minors of $C$ are strictly positive and $\text{det}(C) = 0$. 
By the general theory in \cite{kac}, $C$ is {\it symmetrizable}, that is there is a diagonal matrix with rational coefficients $D = \text{diag}(r_0,\cdots,r_n)$ such that $DC$ is symmetric. 
The {\it quantum affine algebra} $\U_q(\Glie)$ is defined by generators $k_i^{\pm 1}$, $x_i^{\pm}$ ($0\leq i\leq n$) and 
relations
$$k_ik_j=k_jk_i\text{ , } k_ix_j^{\pm}=q^{\left(\pm r_i C_{i,j}\right)}x_j^{\pm}k_i\text{ , }[x_i^+,x_j^-]=\delta_{i,j}\frac{k_i-k_i^{-1}}{q^{r_i}-q^{-r_i}},$$
$$\underset{r=0\cdots 1-C_{i,j}}{\sum}(-1)^r(x_i^{\pm})^{\left(1-C_{i,j}-r\right)}x_j^{\pm}(x_i^{\pm})^{(r)}=0 \text{ (for $i\neq j$)}
,$$
where we denote $\left(x_i^{\pm}\right)^{(r)} = \left(x_i^{\pm}\right)^r/[r]_{q^{r_i}}!$ for $r\geq 0$. We use the standard $q$-factorial notation $[r]_q ! = [r]_q [r-1]_q \cdots [1]_q = (q^r - q^{-r})(q^{r-1} - q^{1 - r})\cdots (q - q^{-1}) (q - q^{-1})^{-r}$. The $x_i^\pm$, $k_i^{\pm 1}$ are called {\it Chevalley generators}.

We use the coproduct $\Delta : \U_q(\Glie)\rightarrow \U_q(\Glie)\otimes \U_q(\Glie)$ defined for $0\leq i\leq n$ by
$$\Delta(k_i)=k_i\otimes k_i\text{ , }\Delta(x_i^+)=x_i^+\otimes 1 + k_i\otimes x_i^+\text{ , }\Delta(x_i^-)=x_i^-\otimes k_i^{-1} + 1\otimes x_i^-.$$ 
This is the same choice as in \cite{da, Chari2, Fre2}\footnote{In \cite{mk} another coproduct is used. We recover the coproduct used in the present paper by taking the opposite coproduct and changing $q$ to $q^{-1}$.}.

\subsection{}\label{rappel}

The indecomposable affine Cartan matrices are classified \cite{kac} into two main classes,
{\it twisted} types and {\it untwisted} types. The latest includes {\it simply-laced} types
and {\it untwisted non simply-laced} types. The type of $C$ is denoted by $X$. We use the numbering of nodes as in \cite{kac} if $X\neq A_{2n}^{(2)}$, and we use the reversed numbering if $X = A_{2n}^{(2)}$.  

 We set $\mu_i = 1$ for $0\leq i\leq n$, except when $(X,i) = (A_{2n}^{(2)}, n)$ where we set $\mu_n = 2$. Without loss of generality, we can choose the $r_i$ so that $\mu_i r_i\in \NN^*$ for any $i$ and $\left(\mu_0 r_0\wedge\cdots\wedge \mu_n r_n\right) = 1$ (there is a unique such choice).

 Let $I = \{1,\cdots, n\}$ and let $\overline{\Glie}$ be the finite-dimensional simple Lie algebra of Cartan matrix $(C_{i,j})_{i,j\in I}$. We denote respectively by $\omega_i$, $\alpha_i$, $\alpha_i^\vee$ ($i\in I$) the fundamental weights, the simple roots and the simple coroots of $\overline{\Glie}$.
We use the standard partial ordering $\leq$ on the weight lattice $P$ of $\overline{\Glie}$. 
The subalgebra of $\U_q(\Glie)$ generated by the $x_i^\pm, k_i^{\pm 1}$ ($i\in I$) is isomorphic to the quantum group of finite type $\U_q\left(\overline{\Glie}\right)$ if $X\neq A_{2n}^{(2)}$, and to $\U_{q^{\frac{1}{2}}}\left(\overline{\Glie}\right)$ if $X = A_{2n}^{(2)}$.  By abuse of notation this algebra will be denoted by $\U_q\left(\overline{\Glie}\right)$.

$\U_q(\Glie)$ has another set of generators, called {\it Drinfeld generators}, denoted by 
$$x_{i,m}^\pm\text{ , }k_i^{\pm 1}\text{ , }h_{i,r}\text{ , }c^{\pm 1/2}\text{ for $i\in I$, $m\in \ZZ$, $r\in\ZZ \setminus \{0\}$,}$$ 
and defined from the Chevalley generators by using the action of Lusztig automorphisms of $\U_q(\Glie)$ (in \cite{bec} for the untwisted types and in \cite{da} for the twisted types).
We have $x_i^\pm = x_{i,0}^\pm$ for $i\in I$. For the untwisted types, a complete set of relations have been proved for the Drinfeld generators \cite{bec, bcp}. For the twisted types, only a partial set of relations have been established (at the time this paper is written), but they are sufficient to study finite-dimensional representations (see the discussion and references in \cite{h8}).
In particular for all types the multiplication defines a surjective linear morphism
\begin{equation}\label{trian}\U_q^-(\Glie)\otimes \U_q(\Hlie)\otimes \U_q^+(\Glie)\rightarrow \U_q(\Glie)\end{equation}
where $\U_q^\pm(\Glie)$ is the subalgebra generated by the $x_{i,m}^\pm$ ($i\in I$, $m\in\ZZ$) and $\U_q(\Hlie)$ is the subalgebra generated by the $k_i^{\pm 1}$, the $h_{i,r}$ and $c^{\pm 1/2}$ ($i\in I$, $r\in\ZZ\setminus\{0\}$).

For $i\in I$, the action of $k_i$ on any object of $\mathcal{F}$ is diagonalizable with eigenvalues in $\pm q^{r_i\ZZ}$.
Without loss of generality, we can assume that $\mathcal{F}$ is the category of {\it type 1} finite-dimensional representations (see \cite{Cha2}), i.e. we assume that for any object of $\mathcal{F}$, the eigenvalues of $k_i$ are in $q^{r_i \ZZ}$ for $i\in I$.

For the untwisted types, the simple objects of $\mathcal{F}$ have been classified by Chari-Pressley \cite{cp, Cha2} by using the Drinfeld generators. For the twisted types, the proof is given in \cite{Cha5, h8}. In both cases the simple objects are parameterized by $n$-tuples of polynomials $(P_i(u))_{i\in I}$ satisfying $P_i(0) = 1$ (they
are called {\it Drinfeld polynomials}).

The action of $c^{\pm 1/2}$ on any object $V$ of $\mathcal{F}$ is the identity, and so the action of
the $h_{i,r}$ commute. Since they also commute with the $k_i$, $V$ can be decomposed 
into generalized eigenspaces $V_m$ for the action of all the $h_{i,r}$ and all the $k_i$ :
$$V = \bigoplus_{m\in\mathcal{M}} V_m.$$  

The $V_m$ are called {\it $l$-weight spaces}. By the Frenkel-Reshetikhin $q$-character theory \cite{Fre}, the eigenvalues can be {\it encoded} by {\it monomials} $m$ in formal variables $Y_{i,a}^{\pm 1}$ ($i\in I, a\in\CC^*$). The construction is extended to twisted types in \cite{h8}.
$\mathcal{M}$ is the set of such monomials (also called {\it $l$-weights}). The {\it $q$-character morphism} is an injective ring morphism
$$\chi_q : \text{Rep}(\U_q(\Glie)) \rightarrow \Yim = \ZZ\left[Y_{i,a}^{\pm 1}\right]_{i\in I, a\in\CC^*},$$
$$\chi_q(V) = \sum_{m\in \mathcal{M}} \text{dim}(V_m) m.$$

For the twisted types there is a modification of the theory and we consider two kinds of variables in \cite{h8}. For homogeneity of notations, the $Y_{i,a}$ in the present paper are the $Z_{i,a}$ of \cite{h8} (we do not use in this paper the variables denoted by $Y_{i,a}$ in \cite{h8}, so there is no possible confusion).

\begin{rem}\label{sumprod} For any $i\in I$, $r\in\ZZ\setminus\{0\}$, $m,m'\in\mathcal{M}$, the eigenvalue of $h_{i,r}$ associated to $mm'$ is the sum of the eigenvalues of $h_{i,r}$ associated respectively to $m$ and $m'$ \cite{Fre, h8}.
\end{rem}

If $V_m\neq \{0\}$ we say that $m$ is an {\it $l$-weight of $V$}.
A vector $v$ belonging to an $l$-weight space $V_m$ is called an {\it $l$-weight vector}. We denote $M(v) = m$ the {\it $l$-weight of $v$}.
A {\it highest $l$-weight vector} is an $l$-weight vector $v$ satisfying $x_{i,p}^+ v = 0$ for any $i\in I$, $p\in\ZZ$.

For $\omega\in P$, the {\it weight space} $V_\omega$ is the set of {\it weight vectors} of weight $\omega$ for $\U_q(\overline{\Glie})$, i.e. of vectors $v\in V$ satisfying $k_i v = q^{\left(r_i \omega(\alpha_i^\vee)\right)} v$ for any $i\in I$.
We have a decomposition $V = \bigoplus_{\omega \in P} V_\omega$. The decomposition in $l$-weight spaces is finer than the decomposition in weight spaces. Indeed, if $v\in V_m$, then $v$ is a weight vector of weight 
$$\omega(m) = \sum_{i\in I, a\in\CC^*} u_{i,a}(m) \mu_i\omega_i\in P,$$ 
where we denote $m = \prod_{i\in I, a\in\CC^*}Y_{i,a}^{u_{i,a}(m)}$. For $v\in V_m$, we set $\omega(v) = \omega(m)$.

A monomial $m\in\mathcal{M}$ is said to be {\it dominant} if $u_{i,a}(m)\geq 0$ for any $i\in I, a\in\CC^*$. For $V$ a simple object in $\mathcal{F}$, let $M(V)$ be the {\it highest weight monomial} of $\chi_q(V)$ (that is $\omega(M(V))$ is maximal). $M(V)$ is dominant and characterizes the isomorphism class of $V$ (it is equivalent to the data of the Drinfeld polynomials). Hence to a dominant monomial $M$ is associated a simple representation $L(M)$. 

\begin{ex}\label{basic} Let us recall the following standard example \cite{j, cp} which we shall use in the following. Let $a\in\CC^*$ and let $L_a = \CC v_a^+ \oplus \CC v_a^-$ be the fundamental representation of $\U_q(\hat{sl_2})$ with the action of the Chevalley generators recalled in the following table.
$$\begin{array}{l|l|l|l|l|l|l}
 & x_1^+ & x_1^- & x_0^+ & x_0^- & k_1 & k_0
\\ v_a^+ & 0 & v_a^- & a v_a^- & 0 & q v_a^+ & q^{-1} v_a^+
\\ v_a^- & v_a^+ & 0 & 0 & a^{-1} v_a^+ & q^{-1} v_a^- & q v_a^- 
\end{array}$$
We have $h_{1,1} = q^{-2} x_1^+x_0^+ - x_0^+x_1^+$ \cite{bec}. The eigenvalue of $h_{1,1}$ corresponding to $m\in\mathcal{M}$ is $\sum_{a\in\CC^*}u_{1,a}(m) a$ \cite{Fre}. We get $h_{1,1}.v_a^+ = aq^{-2}v_a^+$, $k_1.v_a^+ = qv_a^+$, so $M(v_a^+) = Y_{1,aq^{-2}}$. In the same way $M(v_a^-) = Y_{1,a}^{-1}$. Hence $\chi_q(L_a) = Y_{1,aq^{-2}} + Y_{1,a}^{-1}$ and $L_a = L\left(Y_{1,aq^{-2}}\right)$.
\end{ex}

Let $i\in I, a\in\CC^*$ and let us define the monomial $A_{i,a}$ analog of a simple root.
For the untwisted cases, we set \cite{Fre}
$$A_{i,a} = Y_{i,aq^{-r_i}}Y_{i,aq^{r_i}}\times \left(\prod_{\left\{j\in I|C_{i,j} = -1\right\}}Y_{j,a}\right)^{-1}$$
$$\times\left(\prod_{\left\{j\in I|C_{i,j} = -2\right\}}Y_{j,aq^{-1}}Y_{j,aq}\right)^{-1}\times\left(\prod_{\left\{j\in I|C_{i,j} = -3\right\}}Y_{j,aq^{-2}}Y_{j,a}Y_{j,aq^2}\right)^{-1}.
$$

For the twisted cases, let $r$ be the {\it twisting order} of $\Glie$, that is $r = 2$ if $X\neq D_4^{(3)}$ and $r = 3$ if $X = D_4^{(3)}$.
Let $\epsilon$ be a primitive $r$th root of $1$ (for the untwisted cases we set by convention $r = 1$ and $\epsilon = 1$). We now define $A_{i,a}$ as in \cite{h8}. 
\\If $(X,i)\neq (A_{2n}^{(2)}, n)$ and $r_i = 1$, we set
$$A_{i,a} = Y_{i,aq^{-1}}Y_{i,aq}\times\left(\prod_{\left\{j\in I | C_{i,j} < 0\right\}}Y_{j,a^{\left(r_jC_{j,i}\right)}}\right)^{-1}.$$
If $(X,i)\neq (A_{2n}^{(2)}, n)$ and $r_i > 1$, we set
$$A_{i,a} = Y_{i,aq^{-r_i}}Y_{i,aq^{r_i}}\times\left(\prod_{\left\{j\in I | C_{i,j} < 0,r_j = r\right\}}Y_{j,a}\right)^{-1}
\times\left(\prod_{\left\{j\in I | C_{i,j} < 0, r_j = 1\right\}}\left(\prod_{\left\{b\in\CC^*|(b)^r = a\right\}} Y_{j,b}\right)\right)^{-1}.$$
If $(X,i) = (A_{2n}^{(2)}, n)$, we set $
A_{n,a} = \begin{cases}Y_{n,aq^{-1}}Y_{n,aq}Y_{n,-a}^{-1}Y_{n-1,a}^{-1}\text{ if $n > 1$,} 
\\Y_{1,aq^{-1}}Y_{1,aq}Y_{1,-a}^{-1}\text{ if $n = 1$.}\end{cases}$

As mentioned in the introduction, the statement of Theorem \ref{factg} is proved when $\Glie = \hat{sl}_2$ by Chari-Pressley. Let us explain the proof in this case. For $k\geq 1$, $a\in\CC^*$, let $W_{k,a} = L(Y_{1,a} Y_{1,aq^2}\cdots Y_{1,aq^{2(k-1)}})$ (this is called a {\it Kirillov-Reshetikhin module}). A {\it $q$-segment} is a subset of $\CC^*$ of the form $\left\{z, zq^2,\cdots, zq^{2K}\right\} = \left[z, zq^{2K}\right]$ where $z\in\CC^*$ and $K\in \NN$. We use the notation with $\NN$ containing $0$.
We say that $W_{k,a}$ and $W_{l,b}$ are in {\it special position} if $\left[a,aq^{2(k-1)}\right]\cup \left[b,bq^{2(l-1)}\right]$ is a $q$-segment which contains properly $\left[a,aq^{2(k-1)}\right]$ and $\left[b,bq^{2(l-1)}\right]$. Otherwise they are said to be in {\it general position}. 
\begin{thm}\label{fact}\cite{cp, Cha2} The tensor product $W_{k_1,a_1}\otimes \cdots \otimes W_{k_L,a_L}$ is simple if and only if $W_{k_i,a_i}$ and $W_{k_j,a_j}$ are in general position for any $1\leq i < j\leq L$. Any simple object of $\mathcal{F}$ can be factorized in this form.
\end{thm}

 This is an explicit factorization in prime representations. Now the statement of Theorem \ref{factg} for $\Glie = \hat{sl}_2$ follows immediately. Notice that Theorem \ref{fact} also implies that for $\Glie = \hat{sl}_2$, any simple object of $\mathcal{F}$ is real, as mentioned in the introduction.

 For $\Glie \neq \hat{sl}_2$, such a nice description of the factorization is not known. That is why the proof in the present paper does not involve such kind of factorizations. 

\section{Tensor products of $l$-weight vectors}\label{deux}\label{tplwv}

In this section we prove a general result for tensor products of $l$-weight vectors (Proposition \ref{prodlweight}).

 $\U_q(\Glie)$ has a natural grading by the
root lattice $Q = \sum_{i\in I}\ZZ\alpha_i$ of $\overline{\Glie}$ defined by
$$\deg\left(x^\pm_{i,m}\right) = \pm\a_i\text{ , }
\deg\left(h_{i,r}\right)=\deg\left(k_i^\pm\right) = \deg\left(c^{\pm 1/2}\right) = 0.$$

Let $\tilde{\U}_q^+(\Glie)$ (resp. $\tilde{\U}_q^-(\Glie)$) be the subalgebra of $\U_q(\Glie)$
consisting of elements of positive 
(resp. negative) $Q$-degree. These subalgebras should not be confused with the subalgebras $\U_q^\pm(\Glie)$ previously defined in terms of Drinfeld generators. Let 
$$X^+ = \sum_{j\in I, m\in\ZZ}\CC x_{j,m}^+\subset \tilde{\U}_q^+(\Glie).$$

\begin{thm}\label{apco} Let $i\in I$, $r > 0$, $m\in\ZZ$. We have
\begin{equation}\label{h}\Delta\left(h_{i,r}\right) \in h_{i,r}\otimes 1 + 1 \otimes h_{i,r} +  \tilde{\U}_q^-(\Glie) \otimes \tilde{\U}_q^+(\Glie),\end{equation}
\begin{equation}\label{x}\Delta\left(x_{i,m}^+\right) \in x_{i,m}^+\otimes 1 + \U_q(\Glie)\otimes \left(\U_q(\Glie) X^+\right).\end{equation}
\end{thm}

 For the untwisted types the proof can be found in \cite[Proposition 7.1]{da}. For $X = A_2^{(2)}$ see \cite{Cha5} and for the general twisted types see \cite[Proposition 7.1.2]{da2}, \cite[Proposition 7.1.5]{da2}, \cite[Theorem 2.2]{jm}. 

Let $\U_q(\Hlie)^+$ be the subalgebra of $\U_q(\Hlie)$ generated by the $k_i^{\pm 1}$ and the $h_{i,r}$ ($i\in I, r > 0$). The $q$-character and the decomposition in $l$-weight spaces of a representation in $\mathcal{F}$ is completely determined by the action of $\U_q(\Hlie)^+$ \cite{Fre, h8}. Therefore one can define the $q$-character $\chi_q(W)$ of a $\U_q(\Hlie)^+$-submodule $W$ of an object in $\mathcal{F}$.

\begin{prop}\label{prodlweight}
Let $V_1, V_2$ in $\mathcal{F}$ and consider an $l$-weight vector
$$w = \left(\sum_{\alpha} w_\alpha\otimes v_\alpha\right) + \left(\sum_\beta w_\beta'\otimes v_\beta'\right)\in V_1\otimes V_2$$ 
satisfying the following conditions.

(i) The $v_\alpha$ are $l$-weight vectors and the $v_\beta'$ are weight vectors.

(ii) For any $\beta$, there is an $\alpha$ satisfying $\omega(v_\beta') > \omega(v_\alpha)$.

(iii) For $\omega\in\{\omega(v_\alpha)\}_\alpha$, we have $\sum_{\left\{\alpha|\omega(v_{\alpha}) = \omega\right\}} w_{\alpha}\otimes v_{\alpha}\neq 0$. 

\noindent Then $M(w)$ is the product of one $M(v_\alpha)$ by an $l$-weight of $V_1$.
\end{prop}

\demo Consider 
$$V = V_1\otimes\left(\bigoplus_{\left\{\omega\in P|\exists \alpha,\omega(v_\alpha)\leq\omega\right\}}\left(V_2\right)_{\omega}\right) \supset \tilde{V} = V_1\otimes \left(\bigoplus_{\left\{\omega\in P|\exists \alpha, \omega(v_\alpha)<\omega\right\}}\left(V_2\right)_{\omega}\right).$$
By Formula (\ref{h}), $V$ and $\tilde{V}$ are sub $\U_q(\Hlie)^+$-module of $V_1\otimes V_2$.
By condition (ii), $w\in V$ and the image of $w$ in $V/\tilde{V}$ is equal to the image $u$ of $\sum_\alpha w_\alpha\otimes v_\alpha$ in $V/\tilde{V}$. By considering $\alpha_0$ such that $\omega(v_{\alpha_0})$ is minimal, it follows from condition (iii) that $\sum_\alpha w_\alpha\otimes v_\alpha\notin \tilde{V}$. Hence $u\neq 0$ and $M(w) = M(u)$.
Again by Formula (\ref{h}), the action of $h_{i,r}$ on $V/\tilde{V}$ is the action of $h_{i,r}\otimes 1 + 1\otimes h_{i,r}$ for any $i\in I$, $r > 0$. These operators commute with all the operators $h_{i,r}\otimes 1$, $1\otimes h_{i,r}$ ($i\in I, r > 0$), which also commute all together. 

Consider $W = \U_q(\Hlie)^+.u \subset (V/\tilde{V})_{\omega(u)}$. As $W$ is finite-dimensional, there is $u'$ in $W$ which is a common eigenvector for the three families $\left(h_{i,r}\otimes 1 + 1\otimes h_{i,r}\right)_{i\in I, r >0}$, $\left(h_{i,r}\otimes 1\right)_{i\in I, r > 0}$, $\left(1\otimes h_{i,r}\right)_{i\in I, r > 0}$.
We get immediately that the eigenvalues of $u'$ for the first two families are encoded respectively
by $M(u)$ and by an $l$-weight $m$ of $V_1$.
By condition (i), each $w_\alpha\otimes v_\alpha$ is a common generalized eigenvector for $\left(1\otimes h_{i,r}\right)_{i\in I , r > 0}$. Hence $W = \sum_\alpha W_\alpha$ where $W_\alpha$ is the space of common generalized eigenvectors for $\left(1\otimes h_{i,r}\right)_{i\in I, r > 0}$ in $W$ with eigenvalues encoded by $M(v_\alpha)$. So there is an $\alpha$ such that $u'\in W_\alpha$. By Remark \ref{sumprod}, we get $M(v_{\alpha}) m = M(u)$ and so $M(v_{\alpha}) m = M(v)$.
\qed

\begin{ex}\label{exun} We use notations and computations as in Example \ref{basic}.
Let $a\neq b\in\CC^*$ and consider $L_a\otimes L_b$ (this is a generalization of \cite[Example 8.4]{hl}). 
We set $w_a^\pm = v_a^\pm$.
We have
$$\chi_q(L_a\otimes L_b) = Y_{1,aq^{-2}}Y_{1,bq^{-2}} + Y_{1,aq^{-2}}Y_{1,b}^{-1} + Y_{1,a}^{-1} Y_{1,bq^{-2}} + Y_{1,a}^{-1}Y_{1,b}^{-1}.$$
We shall find an $l$-weight vector $w$ in $L_a\otimes L_b$ of $l$-weight $Y_{1,aq^{-1}}Y_{1,b}^{-1}$ illustrating Proposition \ref{prodlweight}. First let us give a generator of each $l$-weight space (they are of dimension $1$ as $a\neq b$). $w_a^+\otimes v_b^+$ (resp. $w_a^-\otimes v_b^-$) is of $l$-weight $Y_{1,aq^{-2}}Y_{1,bq^{-2}}$ (resp. $Y_{1,a}^{-1}Y_{1,b}^{-1}$). Let us look at the weight space of weight $0$. The matrix of $h_{1,1}$ on the basis $(w_a^-\otimes v_b^+, w_a^+\otimes v_b^-)$ is
$\begin{pmatrix}
q^{-2}b - a & a(-q + q^{-3})
\\ 0 & (q^{-2} a - b)
\end{pmatrix}$.
Thus, $w_a^-\otimes v_b^+$ has $l$-weight $Y_{1,bq^{-2}}Y_{1,a}^{-1}$ and
$$w = (b-a)(w_a^+\otimes v_b^-) + a (q - q^{-1}) (w_a^-\otimes v_b^+)$$ 
has $l$-weight $Y_{1,aq^{-2}}Y_{1,b}^{-1}$.
Then $w$ satisfies the conditions of Proposition \ref{prodlweight} with a unique $\alpha$, a unique $\beta$,
$w_\alpha = (b-a)w_a^+$, $v_\alpha = v_b^-$, $w_\beta' = a (q - q^{-1}) w_a^-$ and $v_\beta' =  v_b^+$. The $l$-weight of $w$ is equal to the product $M(w) = Y_{1,aq^{-2}}Y_{1,b}^{-1}$ of $M(v_\alpha) = Y_{1,b}^{-1}$ 
and of $M(w_\alpha)=Y_{1,aq^{-2}}$ which is an $l$-weight of $L_a$.
\end{ex}

\begin{rem} If we replace $v_\alpha$, $v_\beta'$, $V_1$ respectively by $w_\alpha$, $w_\beta'$, $V_2$ in conditions (i), (ii) and in the conclusion, the result does not hold.
For instance $w$ in Example \ref{exun} would satisfy these hypothesis with
$w_\alpha = a (q - q^{-1}) w_a^-$, $v_\alpha =  v_b^+$, 
$w_\beta' = (b-a)w_a^+$, $v_\beta' = v_b^-$. But $M(w)$ is not the product of $M(w_\alpha) = Y_{1,a}^{-1}$ 
by an $l$-weight of $L_b$. 
The reason is that Formula (\ref{h}) also holds for $r < 0$ in the same form, with
a remaining term in $\tilde{\U}_q^-(\Glie) \otimes \tilde{\U}_q^+(\Glie)$ and not in $\tilde{\U}_q^+(\Glie)\otimes \tilde{\U}_q^-(\Glie)$ (this is clear from the relation between the involution $\Omega$ and the coproduct in \cite[Remark 6,(5)]{da}). 
\end{rem}

\section{Reduction and involution}\label{bis}\label{reduc} 

In this section we reduce the proof of Theorem \ref{factg}.

\subsection{} In this subsection we shall review general results which are already known for the untwisted types. For completeness we also give the proofs for the twisted types.

Let $i\in I$. If $r = r_i > 1$ we set $d_i = r_i$. We set $d_i = 1$ otherwise. So for the twisted types we have $d_i = \mu_ir_i$, and for the untwisted types we have $d_i = 1$. We define the fundamental representation $V_i(a) = L\left(Y_{i,a^{d_i}}\right)$ for $i\in I$, $a\in\CC^*$.

If $\Glie$ is twisted, let $\tilde{\Glie}$ be the simply-laced affine Lie algebra associated to $\Glie$ \cite{kac}. Let $\tilde{I}$ be the set of nodes of the Dynkin diagram of the underlying finite-dimensional Lie algebra, with its twisting $\sigma : \tilde{I}\rightarrow \tilde{I}$ and the projection $\tilde{I} \rightarrow I$. We choose a connected set of representatives so that we get $I\subset \tilde{I}$ by identification. 
To avoid confusion, the fundamental representations of $\U_q(\tilde{\Glie})$ are denoted by $\tilde{V}_i(a)$,
the $q$-character morphism of $\U_q(\tilde{\Glie})$ by $\tilde{\chi}_q$, and the corresponding variables by $\tilde{Y}_{i,a}^{\pm 1}$. Consider the ring morphism
$$\pi : \ZZ\left[\tilde{Y}_{i,a}^{\pm 1}\right]_{i\in \tilde{I}, a\in\CC^*}\rightarrow \ZZ\left[Y_{i,a}^{\pm 1}\right]_{i\in I, a\in\CC^*},$$
$$\pi\left(\tilde{Y}_{\sigma^p(i),a}\right) = Y_{i, \left(\epsilon^p a\right)^{d_i}}\text{ for $i\in I, a\in\CC^*, p\in\ZZ$.}$$

\begin{prop}\label{ffund}\cite[Theorem 4.15]{h8} Let $i\in I$, $a\in\CC^*$. We have
$$\chi_q\left(V_i(a)\right) = \pi\left(\tilde{\chi}_q\left(\tilde{V}_i(a)\right)\right).$$
\end{prop} 

\begin{lem}\label{fund} Let $i\in I$, $a\in\CC^*$. We have
$$\chi_q(V_i(a)) \in Y_{i,a^{d_i}} + Y_{i,a^{d_i}}A_{i,\left(a^{d_i} q^{\mu_i r_i}\right)}^{-1}\ZZ\left[A_{j,\left(a\epsilon^k q^r\right)^{d_j}}^{-1}\right]_{j\in I, k\in\ZZ,r>0}.$$ 
\end{lem}

\demo For the untwisted types the proof can be found in the proof of \cite[Lemma 6.5]{Fre2}. For the twisted types, the result follows from Proposition \ref{ffund}.
\qed

For $m, m'\in\mathcal{M}$, we set $m\leq m'$ if $m\in m'\ZZ\left[A_{i,a}^{-1}\right]_{i\in I,a\in\CC^*}$. This defines a partial ordering on $\mathcal{M}$ as the $A_{i,a}^{-1}$ are algebraically independent \cite{Fre}. 

\begin{prop}\label{lower} Let $m\in\mathcal{M}$ dominant. We have
$$\chi_q(L(m)) \in m + \sum_{m' < m} \ZZ m'.$$
\end{prop}

\demo For the untwisted types the result is proved in \cite[Theorem 4.1]{Fre2}. In general it is a direct consequence of Lemma \ref{fund} since a simple module is a subquotient of a tensor product of fundamental representations \cite{Cha2, Cha5}.
\qed

For $a\in\CC^*$, consider the ring (which depends only on the class of $a$ in $\CC^*/\left(q^{\ZZ}\epsilon^\ZZ\right)$)
$$\Yim_a = \ZZ\left[Y_{i,\left(aq^l\epsilon^k\right)^{d_i}}^{\pm 1}\right]_{i\in I, l,k\in\ZZ}.$$
\begin{rem}\label{exya} For instance, we have $A_{i,\left(aq^l\epsilon^k\right)^{d_i}}\in\Yim_a$ for any $i\in I$, $l,k\in\ZZ$. Consequently, $\chi_q(V_i(a))\in\Yim_a$ for any $i\in I$, since $\mu_ir_i\in d_i\ZZ$.
\end{rem}

\begin{defi} $\mathcal{C}$ is the full subcategory of objects in $\mathcal{F}$ whose Jordan-H\"older composition series involve simple representations $V$ satisfying $M(V)\in \Yim_1$. 
\end{defi}

\noindent For instance, the representation $V_i\left(q^l\epsilon^k\right)$ is an object of $\mathcal{C}$ for any $i\in I$, $l,k\in\ZZ$. 
\begin{thm} Let $m\in\mathcal{M}$ dominant. Then there is a unique factorization
$$L(m) \cong \bigotimes_{a\in\left(\CC^*/\left(q^{\ZZ}\epsilon^{\ZZ}\right)\right)} L(m_a)\text{ where }\chi_q(L(m_a))\in\Yim_a.$$\end{thm}
This is a well-known result. The irreducibility of the tensor product follows for example from the criterion in \cite{Chari2} (in other words, it can be proved as in Corollary \ref{map} below that the tensor product and its dual are cyclic). Note that $m_a\in\Yim_a$ implies $\chi_q(L(m_a))\in\Yim_a$ since this holds for fundamental representations by Remark \ref{exya}. 

As a consequence, we can assume in the proof of Theorem 
\ref{factg} that $S_1,\cdots, S_N$ are objects of $\mathcal{C}$. Hence, in the rest of this paper we work in the category $\mathcal{C}$.

Given $a\in\CC^*$, there exists a unique algebra automorphism $\tau_a:\U_q(\Glie)\rightarrow \U_q(\Glie)$ which is defined \cite{cp} on the Drinfeld generators by 
$$\tau_a\left(x^\pm_{i,m}\right)=a^{\pm m} x^\pm_{i,m},\ \
\tau_a\left(h_{i,r}\right)=a^r h_{i,r},\ \ \tau_a\left(k_i^{\pm
1}\right)=k_i^{\pm 1},\ \ \tau_a\left(c^{\pm\frac{1}{2}}\right) = c^{\pm \frac{1}{2}}.$$
The definition in \cite{cp} is given for the untwisted types, but it holds for the twisted types as well. For $V$ in $\mathcal{F}$, let $\tau_a^*(V)$ be the object in $\mathcal{F}$ obtained by pulling back $V$ via the automorphism $\tau_a$.  
The following is proved in \cite{Fre2} for the untwisted types (the highest weight term is computed in \cite{Chari2}). The proof is the same for the twisted types.

\begin{prop}\label{shift} $\chi_q(\tau_a^*(V))$ is obtained from $\chi_q(V)$
by changing each $Y_{i,b}^{\pm 1}$ to $Y_{i,\left(a^{d_i}b\right)}^{\pm 1}$.
\end{prop}

Let $\sigma$ be the involution of the algebra $\Yim$ defined by $Y_{i,a}^{\pm 1}\mapsto Y_{i,a^{-1}}^{\mp 1}$. 

\noindent If $\Glie$ is untwisted, let $w_0$ be the longest element in the Weyl group of $\overline{\Glie}$ and $i\mapsto \overline{i}$ the unique bijection of $I$ satisfying $w_0(\alpha_i) = -\alpha_{\overline{i}}$. Let $h^{\vee}$ be the dual Coxeter number of $\overline{\Glie}$ and $r^{\vee}$ the maximal number of edges connecting two vertices of its Dynkin diagram.

\noindent If $\Glie$ is twisted, we use the same definition of $r^{\vee}, h^{\vee}$, but for $\overline{(\tilde{\Glie})}$. We set $\overline{i} = i$. 

\noindent Let $i\in I$, $a\in \CC^*$. We set $\overset{\circ}{Y}_{i,a} = Y_{i,-a}$ if $X = A_n^{(2)}$ and $r_i\leq 1$. We set $\overset{\circ}{Y}_{i,a} = Y_{\overline{i},a}$ otherwise.

The proof of the following duality result can be found in \cite{h4} for the untwisted types (the highest weight term is computed in \cite{cp5}).

\begin{prop}\label{invmon} For $m\in \Yim_1$ a dominant monomial, we have
$$\sigma\left(\chi_q(L(m))\right) = \chi_q\left(L\left(M\right)\right)\text{ where }M = \prod_{i\in I, a\in (\epsilon^\ZZ q^{\ZZ})^{d_i}}\left(\overset{\circ}{Y}_{i,\left(a^{-1} q^{-\left(d_i r^{\vee}h^{\vee}\right)}\right)}\right)^{u_{i,a}(m)}.$$
\end{prop}

\demo It suffices to prove the statement for the twisted types. By using the arguments of \cite[Section 4.2]{h4}, it suffices to compute the lowest weight monomial of $\chi_q(L(m))$. By using the same argument as in \cite[Corollary 6.9]{Fre2}, it suffices to do it for fundamental representations. Hence the result follows from Proposition \ref{ffund}.
\qed

\subsection{}\label{inddepart}

We define a sequence of subcategories of $\mathcal{C}$ in the spirit of the categories in \cite{hl} (but the subcategories that we consider here are different from the categories in \cite{hl}).

\begin{defi} Let $\ell\geq 0$. $\mathcal{C}_\ell$ is the full subcategory of $\mathcal{C}$ whose Jordan-H\"older composition series involves simple representations $V$ satisfying 
$$M(V)\in \ZZ\left[Y_{i,\left(q^l\epsilon^k\right)^{d_i}}\right]_{i\in I,0\leq l\leq \ell,k\in\ZZ}.$$
\end{defi}

For instance the representation $V_i(q^l\epsilon^k)$ is an object of $\mathcal{C}$ for $i\in I$, $k\in\ZZ$, $0\leq l\leq \ell$.

A monomial $m\neq 1$ in $\Yim_1$ is said to be right-negative (resp. left-negative) \cite{Fre2} 
if the factors $Y_{i,\left(q^r\epsilon^k\right)^{d_i}}$ appearing in $m$, for which $r$ is maximal (resp. minimal), have negative powers. The $A_{i,a}^{-1}$ are right-negative and a product of right-negative monomials is right-negative. A right-negative monomial is not dominant \cite{Fre2}.

An analog of the following result was proved in \cite{hl} for simply-laced types.

\begin{lem} $\mathcal{C}_\ell$ is stable under tensor product and $\mathcal{C}_\ell$ inherits a structure of a tensor category.\end{lem}

\demo Let $L(m), L(m')$ be simple objects of $\mathcal{C}_\ell$ and $L(m'')$ be a simple constituent of $L(m)\otimes L(m')$. Then from Lemma \ref{fund}, $m''$ is of the form $mm' A$ where $A$ is a product of $A_{i,\left(\epsilon^k q^r\right)^{d_i}}^{-1}$ with $i\in I$, $k\in\ZZ$, $r > 0$. Suppose that one factor $A_{i,\left(\epsilon^{kd_i} q^R\right)}^{-1}$ occurs in $A$ with $R >  d_i\ell - \mu_i r_i$. Then $mm'A_{i,\left(\epsilon^{kd_i} q^R\right)}^{-1}$ is right-negative, so $mm'A$ is right-negative as a product of right-negative monomials. Contradiction as $m''$ is dominant. Hence the $A_{i,\left(\epsilon^{d_ik}q^r \right) }^{-1}$ occurring in $A$ satisfy $0<r\leq d_i\ell - \mu_i r_i$. So $m''\in \ZZ\left[Y_{i,\left(\epsilon^kq^l\right)^{d_i}}\right]_{i\in I, 0\leq l\leq \ell,k\in\ZZ}$ and $L(m'')$ is an object of $\mathcal{C}_\ell$. Hence $L(m)\otimes L(m')$ is an object of $\mathcal{C}_\ell$.
\qed

The statement of Theorem \ref{factg} is clear for $\mathcal{C}_0$ from the following.

\begin{lem}\label{lzero} All simple objects of $\mathcal{C}_0$ are tensor products of fundamental representations in $\mathcal{C}_0$. An arbitrary tensor product of simple objects in $\mathcal{C}_0$ is simple.\end{lem}

\demo From Lemma \ref{fund}, for any $k\in \ZZ$, $i\in I$, we have
\begin{equation}\label{expfdo}\chi_q\left(V_i(\epsilon^k)\right) \in Y_{i,\epsilon^{k d_i}} + Y_{i,\epsilon^{k d_i}}A_{i,\left(\epsilon^{k d_i}q^{\mu_ir_i}\right)}^{-1}\ZZ\left[A_{j,\left(\epsilon^m q^r\right)^{d_j}}^{-1}\right]_{j\in I, m\in\ZZ, r > 0}.\end{equation}
Let $V$ be a tensor product of such fundamental representations in $\mathcal{C}_0$. By Formula (\ref{expfdo}), a monomial occurring in $\chi_q(V)$ not of highest weight is a product of one $Y_{i,\epsilon^{k d_i}}A_{i,\left(\epsilon^{k d_i}q^{\mu_ir_i}\right)}^{-1}$ by some $Y_{j,\epsilon^{m d_j}}$, $A_{j,\left(\epsilon^m q^r\right)^{d_j}}^{-1}$.
So it is right-negative and not dominant. Hence $V$ is simple and the first statement is proved. 
Let $L(m)$ be a simple object in $\mathcal{C}_0$. Then $m$ is a product of some $Y_{i,\epsilon^{k d_i}}$ ($i\in I, k\in\ZZ$). The second statement follows immediately.
\qed

\begin{rem} The category $\mathcal{C}_0$ is not semi-simple. For instance, for $\Glie = \hat{sl}_2$, $L(Y_{1,1})$ has a non-split self-extension, which can be constructed by a direct computation.\end{rem}

As a consequence of Proposition \ref{shift}, for any simple object $V$ in $\mathcal{C}$, there is
$a\in q^\ZZ\epsilon^\ZZ$ and $\ell\geq 0$ such that $\tau_a^* (V)$ is an object in $\mathcal{C}_\ell$.
Hence it suffices to prove the statement of Theorem \ref{factg} for the categories $\mathcal{C}_\ell$.

Let $\ell \geq 0$. Consider the bar involution defined on $\ZZ\left[Y_{i,\left(\epsilon^k q^l\right)^{d_i}}^{\pm 1}\right]_{i\in I, k\in\ZZ, 0\leq l\leq \ell}$ by $\overline{Y_{i,\left(\epsilon^k q^l\right)^{d_i}}} = Y_{i,\left(\epsilon^{-k} q^{\ell - l}\right)^{d_i}}$ for $i\in I, k\in\ZZ, 0\leq l\leq \ell$. 

For a simple $V = L(m)$ we set $\overline{V} = L(\overline{m})$. This defines a bar involution of the Grothendieck ring of $\mathcal{C}_\ell$. For example, $\overline{V_i(\epsilon^k q^l)} = V_i(\epsilon^{-k} q^{\ell - l})$ for $i\in I$, $k\in\ZZ$, $0\leq l\leq \ell$.

\begin{prop}\label{tourne} For $V$ a simple object in $\mathcal{C}_\ell$ we have $\overline{\chi_q(V)} = \chi_q(\overline{V})$. In particular the bar involution is a ring automorphism of the Grothendieck ring of $\mathcal{C}_\ell$.
\end{prop}

\demo First by using Proposition \ref{invmon} and Proposition \ref{shift}, we get
$$\phi(\chi_q(L(m))) = \chi_q(L(\phi(m)))$$ 
where $\phi$ is the ring isomorphism of $\Yim$ defined by $\phi(Y_{i,a}) = \overset{\circ}{Y}_{i,q^{\ell d_i}a^{-1}}$. Then consider the ring automorphism ${\psi}$ of $\Yim$ defined by $Y_{i,a}\mapsto \overset{\circ}{Y}_{i, a}$. We get immediately
$$\left({\psi}\circ \phi\right)(\chi_q(L(m))) = \chi_q\left(L\left(\left({\psi}\circ \phi\right)(m)\right)\right).$$ 
This is exactly the relation $\overline{\chi_q(L(m))} = \chi_q(L(\overline{m}))$.\qed

\section{Upper and lower $q$-characters}\label{trois}

In this section we introduce the notions of lower and upper $q$-characters that we shall use in the following. We prove several results and formulae about them. Fix $L\in\ZZ$.

\subsection{} For a monomial $m\in\Yim_1$, we denote by $m^{= L}$ (resp. $m^{\leq L}$, $m^{\geq L}$) the product with multiplicities of the factors $Y_{i,\left(\epsilon^k q^l\right)^{d_i}}^{\pm 1}$ occurring in $m$ with $l = L$ (resp. $l\leq L$, $l\geq L)$, $i\in I$, $k\in\ZZ$. 

Consider a dominant monomial $M\in \Yim_1$ and let $V = L(M)$.

\begin{defi}\label{lu} The lower (resp. upper) $q$-character $\chi_{q,\leq L}(V)$ (resp. $\chi_{q,\geq L}(V)$) is the sum with multiplicities of the monomials $m$ occurring in $\chi_q(V)$ satisfying
$$m^{\geq (L+1)} = M^{\geq (L + 1)}\text{ }\left(\text{resp. }m^{\leq (L-1)} = M^{\leq (L - 1)}\right).$$
\end{defi}
\noindent We define $V_{\leq L}, V_{\geq L}\subset V$ as the corresponding respective sums of $l$-weight spaces. 

Let $\mathcal{A}_{\leq L}$ (resp. $\mathcal{A}_{\geq L}$) be the subring of $\Yim$ generated by the 
$A_{i,a}^{-1}$ with $i\in I$ and 
$$a \in \epsilon^{d_i\ZZ}q^{\left(d_i(L - \NN) - \mu_ir_i\right)}\text{ }\left(\text{ resp. }a\in \epsilon^{d_i \ZZ}q^{\left(d_i(L + \NN) + \mu_ir_i\right)}\right).$$

\begin{lem}\label{alternate} 
$\chi_{q,\leq L}(V)$ (resp. $\chi_{q,\geq L}(V)$) is equal to the sum with multiplicities of the monomials $m$ occurring in $\chi_q(V)$ satisfying $mM^{-1}\in \mathcal{A}_{\leq L}$ (resp. $mM^{-1}\in\mathcal{A}_{\geq L}$).
\end{lem}

\demo Let us prove the statement for $\chi_{q,\geq L}(V)$ (the other proof is analog). By Proposition \ref{lower}, we can assume $m\leq M$. First $mM^{-1}\in \mathcal{A}_{\geq L}$ implies $m^{\leq (L-1)} = M^{\leq (L - 1)}$, since for any $i\in I$, $l \geq L$, $k\in\ZZ$, the monomial $A_{i,\left(\left(\epsilon^k q^l\right)^{d_i} q^{\mu_ir_i}\right)}^{-1}$ does not contain any $Y_{j,\left(\epsilon^k q^r\right)^{d_j}}^{\pm 1}$ with $r < L$. To prove the converse, suppose that some $A_{i,\left(\epsilon^k q^l\right)^{d_i}}^{-1}$ with $l < d_iL + \mu_i r_i$ and $k\in\ZZ$ occurs in $mM^{-1}$. Then $mM^{-1}$ is left-negative and there are $j\in I$, $M < L$, $K\in\ZZ$ satisfying $u_{j,\left(q^M\epsilon^{K}\right)^{d_j}}(m M^{-1}) < 0$. Hence $m^{\leq (L-1)} \neq M^{\leq (L-1)}$.
\qed
\begin{rem}\label{plu} As a consequence, for $V,V'$ in $\mathcal{F}$ such that $V\otimes V'$ is simple, we have $\chi_{q,\geq L}(V\otimes V') = \chi_{q,\geq L}(V)\chi_{q,\geq L}(V')$ and $\chi_{q,\leq L}(V\otimes V') = \chi_{q,\leq L}(V)\chi_{q,\leq L}(V')$.\end{rem}

An an illustration, by Lemma \ref{fund} and Lemma \ref{alternate}, for $i\in I$, $k,l\in\ZZ$ we get
$$\chi_{q,\leq l}\left(V_i\left(\epsilon^k q^l\right)\right) = Y_{i,\left(\epsilon^k q^l\right)^{d_i}}\text{ and }\chi_{q,\geq l}\left(V_i\left(\epsilon^k q^l\right)\right) = \chi_q\left(V_i\left(\epsilon^k q^l\right)\right) - Y_{i,\left(\epsilon^k q^l\right)^{d_i}}.$$

\subsection{} A module in $\mathcal{F}$ is said to be {\it cyclic} if it is generated by a highest $l$-weight vector. 

We have the following cyclicity result \cite{Chari2, mk, VaVa3}.

\begin{thm}\label{cyc}\label{thK} Consider $a_1,\cdots, a_R\in \epsilon^\ZZ q^\ZZ$ and $i_1,\cdots, i_R\in I$
such that for $r < p$, we have $a_p a_r^{-1}\in \epsilon^\ZZ q^\NN$. Then the tensor product 
$$V_{i_R}(a_R)\otimes \cdots \otimes V_{i_1}(a_1)$$ 
is cyclic. Moreover there is a unique morphism up to a constant multiple
$$V_{i_R}(a_R)\otimes \cdots \otimes V_{i_1}(a_1) \rightarrow V_{i_1}(a_1)\otimes \cdots \otimes V_{i_R}(a_R),$$
and its image is simple isomorphic to $L\left(Y_{i_1,(a_1)^{d_{i_1}}}\cdots Y_{i_R,(a_R)^{d_{i_R}}}\right)$.
\end{thm}

 Note that the condition in \cite{mk} is that $a_p a_r^{-1}$ has no pole at $q = 0$ when $q$ is an indeterminate. That is why in our context the condition is translated as $a_p a_r^{-1}\in \epsilon^\ZZ q^\NN$. The statement in \cite{mk} involves representations $W(\overline{\omega}_i)\cong V_i(a)$ for some $a\in\CC^*$ computed in \cite[Lemma 4.6]{bn} and \cite[Remark 3.3]{Nad} ($a$ does not depend on $i$ but only on the choice of the isomorphism between Chevalley and Drinfeld realizations).

As a direct consequence of Theorem \ref{cyc}, we get the following.

\begin{cor}\label{map} Let $m,m'\in\Yim_1$ dominant monomials. Assume that $u_{i,a}(m) \neq 0$ implies $u_{j,\left(a(q^{r}\epsilon^k)^{d_i}\right)}(m') = 0$ for any $i,j\in I$, $r > 0$, $k\in\ZZ$, $a\in \left(\epsilon^\ZZ q^\ZZ\right)^{d_i}$. Let $W = L(m)$ and $W' = L(m')$.
Then $W\otimes W'$ is cyclic and there exists a morphism of $\U_q(\Glie)$-modules
$$\mathcal{I}_{W,W'} : W\otimes W' \rightarrow W'\otimes W$$
whose image is simple isomorphic to $L(mm')$.
\end{cor}
This is a well-known result (see for instance \cite{Fre0, ks}). We write the proof for completeness of the paper. The morphism is unique up to a constant multiple. If in addition $W\otimes W'$ is simple, then $W'\otimes W$ is simple as well and $\mathcal{I}_{W,W'}$ is an isomorphism.

\demo From Theorem \ref{cyc}, $W$ (resp. $W'$) is the submodule of a tensor product of fundamental representations $V_1\otimes \cdots \otimes V_R$ (resp. $V_{R+1}\otimes \cdots \otimes V_P$) generated by the tensor product of highest $l$-weight vectors. Consider elements $i_r\in I$, $a_r\in\epsilon^\ZZ q^\ZZ$ satisfying $V_r = V_{i_r}(a_r)$. Then by our assumptions, for $1\leq p < r\leq R$ or $R+1 \leq p < r \leq P$, we have $a_r(a_p)^{-1}\in \epsilon^{\ZZ}q^\NN$. 
Moreover, for $1\leq r\leq R < p\leq P$, we have $a_r(a_p)^{-1}\in \epsilon^\ZZ q^\NN$.
Hence, by Theorem \ref{cyc}, we have surjective morphisms $(V_R\otimes \cdots \otimes V_1)\twoheadrightarrow W$, $(V_P\otimes \cdots \otimes V_{R+1})\twoheadrightarrow W'$ and so a surjective morphism 
$$(V_R\otimes \cdots \otimes V_1)\otimes (V_P\otimes \cdots \otimes V_{R+1})\twoheadrightarrow W\otimes W',$$
where the left-hand module is cyclic. Hence $W\otimes W'$ is cyclic. 
Now, by Theorem \ref{cyc}, for every $1\leq i \leq R < j \leq P$, we have a well-defined morphism $\mathcal{I}_{V_i,V_j} : V_i\otimes V_j\rightarrow V_j\otimes V_i$. So we can consider
$$\mathcal{I} = \left(\mathcal{I}_{V_1,V_P}\circ\cdots \circ \mathcal{I}_{V_R,V_P}\right)\circ \cdots \circ \left(\mathcal{I}_{V_1,V_{R+1}}\circ\cdots \circ \mathcal{I}_{V_R,V_{R+1}}\right) :$$
$$\left(V_1\otimes \cdots \otimes V_R\right)\otimes \left(V_{R+1}\otimes \cdots \otimes V_P\right)\rightarrow \left(V_{R+1}\otimes \cdots \otimes V_P\right)\otimes \left(V_1\otimes \cdots \otimes V_R\right).$$
We use an abuse of notation, as we should have written $\text{Id}\otimes \mathcal{I}_{V_1,V_P}\otimes \text{Id}$. In the following we shall use an analog convention without further comment.
The image of the restriction $\mathcal{I}_{W,W'}$ of $\mathcal{I}$ to $W\otimes W'$ is generated by the tensor product of the highest $l$-weight vectors. Hence it is included in $W'\otimes W$.  By Theorem \ref{cyc} the submodule of $V_{R+1}\otimes \cdots \otimes V_P\otimes V_1\otimes \cdots \otimes V_R$ generated by the tensor product of the highest $l$-weight vectors is simple. Hence the image of $\mathcal{I}_{W,W'}$ is simple.
\qed

\subsection{} We go back to $M\in \Yim_1$ dominant and we turn to studying the surjective morphism
$$\phi = \mathcal{I}_{L\left(M^{\geq L}\right), L\left(M^{\leq (L - 1)}\right)} : L\left(M^{\geq L}\right)\otimes L\left(M^{\leq (L - 1)}\right)\twoheadrightarrow V = L(M).$$
Let $v$ be a highest $l$-weight vector of $L\left(M^{\leq (L - 1)}\right)$.

\begin{prop}\label{useqt} The morphism $\phi$ restricts to a bijection
$$\phi : L\left(M^{\geq L}\right) \otimes v\rightarrow (V)_{\geq L}.$$
\end{prop}

This result generalizes \cite[Lemma 8.5]{hl}. The proof is different because in the general case, the representation
$L\left((M)^{\geq L}\right)$ is not necessarily minuscule (in the sense of \cite{ch}).

\demo First by Formula (\ref{h}), for $i\in I$, $r>0$ and $w\in L\left(M^{\geq L}\right)$, we have
\begin{equation}\label{acth}h_{i,r}(w\otimes v) = \left(h_{i,r}w\right)\otimes v + w \otimes \left(h_{i,r}v\right)\in \left(L\left(M^{\geq L}\right)\otimes v\right).\end{equation} 
Hence $L(M^{\geq L}) \otimes v$ is a $\U_q(\Hlie)^+$-module and by Remark \ref{sumprod}, we get \begin{equation}\label{ftrois}\chi_q\left(L\left(M^{\geq L}\right) \otimes v\right) = \chi_q\left(L\left(M^{\geq L}\right)\right) M^{\leq (L -1)}.\end{equation}
Let us establish
$$\phi^{-1}\left(V_{\geq L}\right) = L\left(M^{\geq L}\right)\otimes v + \text{Ker}(\phi).$$
Clearly $\phi^{-1}\left(V_{\geq L}\right)\supset\text{Ker}(\phi)$. By Formula (\ref{ftrois}), we get $\phi^{-1}\left(V_{\geq L}\right)\supset L\left(M^{\geq L}\right)\otimes v$. So the inclusion $\supset$ is established. Let us prove the other inclusion. 
$\phi^{-1}\left(V_{\geq L}\right)$ is a $\U_q(\Hlie)$-module and so it can be decomposed into $l$-weight spaces 
$$ \phi^{-1}\left(V_{\geq L}\right)= \bigoplus_{m\in\mathcal{M}} \left(\phi^{-1}\left(V_{\geq L}\right)  \right)_m.$$ 
If $m\in\mathcal{M}$ is not an $l$-weight of $V_{\geq L}$, then $\left(\phi^{-1}\left(V_{\geq L}\right)  \right)_m\subset \text{Ker}(\phi)$.
Otherwise, let $w\in \left(\phi^{-1}\left(V_{\geq L}\right)  \right)_m$. If $w\notin L\left((M)^{\geq L}\right)\otimes v$, we can write a decomposition 
$$w = \sum_\alpha w_\alpha\otimes v_\alpha  + w_\beta' \otimes v$$ 
as in Proposition \ref{prodlweight} with all $v_\alpha$ satisfying $\omega(v_\alpha) < \omega(v)$ and only one $v_\beta' = v$. Thus, one of the $M(v_\alpha)$ is a factor of $m$, and so $m$ is not an $l$-weight of $V_{\geq L}$. Contradiction. Hence $\left(\phi^{-1}\left(V_{\geq L}\right)  \right)_m\subset L\left(M^{\geq L}\right)\otimes v$. This concludes the proof of the equality.

Now by Formula (\ref{x}), $L\left((M)^{\geq L}\right)\otimes v$ is stable for the action of the $x_{i,p}^+$, and for $w\in L\left((M)^{\geq L}\right)$, we have 
\begin{equation}\label{actx}x_{i,p}^+(w\otimes v) = \left(x_{i,p}^+ w\right)\otimes v\text{ for any $i\in I$, $p\in\ZZ$.}\end{equation} 
Suppose that there exists a non-zero weight vector $w\otimes v\in \text{Ker}(\phi)\cap \left(L\left(M^{\geq L}\right)\otimes v\right)$. $w\otimes v$ generates a proper submodule of the cyclic module $L(M^{\geq L})\otimes L\left(M^{\leq (L - 1)}\right)$ since $\phi(\U_q(\Glie)(w\otimes v)) = 0$. Let $v'$ be a highest $l$-weight vector of $L(M^{\geq L})$. Since $\omega(w\otimes v) < \omega(v'\otimes v)$, there exists $N\geq 1$ such that there is a decomposition 
$$\omega(w\otimes v) - \omega(v'\otimes v) = - \alpha_{j_1} - \cdots - \alpha_{j_N}\text{ for some $j_1,\cdots,j_N\in I$.}$$ 
Since $L\left(M^{\geq L}\right)\otimes L\left(M^{\leq (L - 1)}\right)$ is cyclic, $v'\otimes v\notin \U_q(\Glie)(w\otimes v)$. Hence for any $i_1,\cdots,i_N\in I$, $p_1,\cdots, p_N\in \ZZ$, we get
\begin{equation}\label{zero}\left(x_{i_1,p_1}^+x_{i_2,p_2}^+\cdots x_{i_N,p_N}^+\right) (w\otimes v) = 0\text{ and }\left(x_{i_1,p_1}^+x_{i_2,p_2}^+\cdots x_{i_N,p_N}^+\right) w = 0.\end{equation}
But $L\left((M)^{\geq L}\right)$ is simple, so there is $g\in\U_q(\Glie)$ satisfying $g w = v'$. By using the surjective map (\ref{trian}), $g$ can be decomposed as a sum of monomials in the Drinfeld generators $g_- h g_+$ where $g_\pm\in\U_q^\pm(\Glie)$ and $h\in \U_q(\Hlie)$. Each term $\left(g_- h g_+\right) w$ is a weight vector and so we can assume that each term satisfies $\omega(g_- h g_+ w) = \omega(v')$. Then each $g_+ w$ is a weight vector satisfying
$\omega(g_+w)\geq \omega(g_- h g_+ w) = \omega(v')$. So each $g_+$ is a product $x_{i_1,p_1}^+x_{i_2,p_2}^+\cdots x_{i_{N'},p_{N'}}^+$ where $N'\geq N$. So $g^+ w = 0$ by Formulae (\ref{zero}). Thus, we have $g w = 0$. Contradiction. Hence we are done since we have established
$$\text{Ker}(\phi)\cap \left(L\left(M^{\geq L}\right)\otimes v\right) = \{0\}.$$\qed

\begin{rem} By Formulae (\ref{acth}), (\ref{actx}), the action of the $x_{i,p}^+$, $h_{i,r}$ on $V_{\geq L}$ can be recovered from their action on $L\left(M^{\geq L}\right)$. This will find other applications in another paper.
\end{rem}

\begin{cor}\label{useqt2} Let $M\in\Yim_1$ be a dominant monomial and $L\in\ZZ$. We have
$$\chi_{q,\geq L}(L(M)) = M^{\leq (L - 1)}\chi_q\left(L\left(M^{\geq L}\right)\right).$$
\end{cor}

\demo In Proposition \ref{useqt}, $\phi$ is an isomorphism of $\U_q(\Hlie)^+$-modules, and so
$$\chi_{q,\geq L}(L(M)) = \chi_q\left(L\left(M^{\geq L}\right) \otimes v\right).$$
We are done by Formula (\ref{ftrois}).\qed

\section{End of the proof of the main theorem}\label{quatre}

First let us mention the proof of the "only if" part of Theorem \ref{factg} (which is trivial). As proved in \cite{Fre}, the injectivity of the $q$-character morphism implies that $\text{Rep}(\U_q(\Glie))$ is commutative (see \cite{h8} for the twisted types). So the irreducibility of $S_1\otimes \cdots \otimes S_N$ is equivalent to the irreducibility of $S_\sigma = S_{\sigma(1)}\otimes \cdots \otimes S_{\sigma(N)}$ for any permutation $\sigma$ of $[1,N]$. Let $i < j$ and $\sigma$ satisfying $\sigma (i) = 1$ and $\sigma (j) = 2$. If $S_i\otimes S_j$ is not simple, we have a proper submodule $V\subset S_i\otimes S_j$ and so a proper submodule $V\otimes S_{\sigma(3)}\otimes \cdots \otimes S_{\sigma(N)}\subset S_\sigma$. Hence $S_1\otimes\cdots \otimes S_N$ is not simple.

Now we turn to the "if" part. We have seen in Section \ref{reduc} that it suffices to prove the statement of Theorem \ref{factg} for the categories $\mathcal{C}_\ell$. We shall proceed by induction on $\ell\geq 0$. For $\ell = 0$ the result has been discussed in Section \ref{inddepart}.

Let $S$ be a simple module in $\mathcal{C}_\ell$ of highest weight monomial $M$. Let $M_- = M^{\leq (\ell - 1)}$ and $M_+ = M^{= \ell }$. Set $S_\pm = L(M_\pm)$. Consider a highest $l$-weight vector $v_\pm$ of $S_\pm$. Recall the surjective morphism of Corollary \ref{map}.
$$\mathcal{I}_{S_+,S_-} : S_+\otimes S_-\twoheadrightarrow S \subset S_-\otimes S_+.$$

\begin{prop}\label{new} Let $S,S'$ simple objects in $\mathcal{C}_\ell$ such that the tensor product $S\otimes S'$ is simple. Then the tensor product $S_-\otimes S'_-$ is simple.
\end{prop}

\demo Let $M = M(S)$ and $M' = M(S')$. As above, we define 
$$M_- = M^{\leq (\ell - 1)}\text{ , }M_+ = M^{= \ell }\text{ , }(M')_- = (M')^{\leq (\ell - 1)}\text{ , }(M')_+ = (M')^{= \ell }.$$ 
The duality of Proposition \ref{tourne} allows to reformulate the problem. Indeed the hypothesis implies that $\overline{S}\otimes \overline{S'}$ is simple, and it suffices to prove that $\overline{S_-}\otimes \overline{S_-'}$ is simple.

From Corollary \ref{useqt2} with $L = 1$, we get 
$$\chi_{q,\geq 1}\left(L\left(\overline{M_+M'_+ M_- M'_-}\right)\right) = \overline{M_+ M'_+}\chi_q\left(L\left(\overline{M_- M'_-}\right)\right).$$
Since $\overline{S}\otimes \overline{S'}$ is simple, this is equivalent to
$$\chi_{q,\geq 1}\left(\overline{S}\otimes \overline{S'}\right) = 
\overline{M_+ M'_+}\chi_q\left(L\left(\overline{M_- M'_-}\right)\right).$$
But by Remark \ref{plu} the left term is equal to $\chi_{q,\geq 1}\left(\overline{S}\right)\chi_{q,\geq 1}\left(\overline{S'}\right)$ which, again by Corollary \ref{useqt2}, is equal to $\overline{M_+ M'_+}\chi_q\left(\overline{S_-}\right)\chi_q\left(\overline{S'_-}\right) = \overline{M_+ M'_+}\chi_q\left(\overline{S_-}\otimes\overline{S'_-}\right)$. This implies 
$$\chi_q\left(\overline{S_-}\otimes \overline{S'_-}\right)
= \chi_q\left(L\left(\overline{M_-M'_-}\right)\right).$$
Hence $\overline{S_-} \otimes \overline{S'_-}$ is simple.
\qed

Now we conclude\footnote{Parts of the final arguments of the present paper were used in the proof of \cite[Theorem 8.1]{hl} for the simply-laced types. But in the context of \cite{hl} the proof is drastically simplified since the $(S_i)_-$ belong to a category equivalent to $\mathcal{C}_0$ and are minuscule.} the proof of Theorem \ref{factg}. In addition to the induction on $\ell$, we start a new induction on $N\geq 2$. For $N = 2$ there is nothing to prove.

For $i=1,\cdots,N$, we define $M_i$, $(M_i)_\pm$, $(S_i)_\pm$, $(u_i)_\pm$ as above.
Consider a pair $(i,j)$ of integers satisfying $1\leq i<j \leq N$.
By our assumptions, $S_i\otimes S_j$ is simple. Hence $(S_i)_-\otimes (S_j)_-$ is simple by Proposition \ref{new}. Besides $(S_j)_+\otimes (S_i)_+$ is a tensor product of fundamental representations belonging to a category equivalent to $\mathcal{C}_0$ by Proposition \ref{shift}. Hence $(S_j)_+\otimes (S_i)_+$ is simple. Now by Corollary \ref{map}, the module
$(S_j)_+\otimes (S_i)_+\otimes (S_i)_-\otimes (S_j)_-$ is cyclic. By Corollary \ref{map}, there exists a surjective morphism 
$$\Psi = \mathcal{I}_{(S_j)_+,(S_j)_-}\mathcal{I}_{(S_j)_+,(S_i)_+}\mathcal{I}_{(S_j)_+,(S_i)_-}\mathcal{I}_{(S_i)_+,(S_i)_-}$$
$$\Psi : (S_j)_+\otimes (S_i)_+\otimes (S_i)_-\otimes (S_j)_- \twoheadrightarrow L(M_i M_j)\cong S_i\otimes S_j.$$
The map $\mathcal{I}_{(S_j)_+,(S_i)_+}\mathcal{I}_{(S_j)_+,(S_i)_-}$ can be rewritten as $\a \otimes {\rm id}_{(S_j)_-}$, where
\[ 
\a : (S_j)_+\otimes (S_i)_-\otimes (S_i)_+ \to (S_i)_-\otimes (S_i)_+\otimes (S_j)_+
\]
restricts to a morphism $\bar{\a} : (S_j)_+\otimes S_i\rightarrow S_i\otimes (S_j)_+$. Now we have 
$$\left(\mathcal{I}_{(S_j)_+,(S_j)_-}\right)^{-1}\left((S_j)_{\geq \ell}\right) = (S_j)_+\otimes (u_j)_-$$
and $\mathcal{I}_{(S_j)_+,(S_j)_-}$ restricts to a bijection from
$(S_j)_+\otimes (u_j)_-$ to $(S_j)_{\geq \ell}$ by Proposition \ref{useqt}.
Since $\Psi$ is surjective, we get
$${\rm Im} (\bar{\a}) \otimes (u_j)_- \supset S_i\otimes (S_j)_+\otimes (u_j)_-.$$ 
Hence $\bar{\a}$ is surjective. 

By the induction hypothesis on $N$, the module $S_1\otimes \cdots \otimes S_{N-1}$ is simple.

Let us prove that $(S_{L})_-\otimes \cdots \otimes (S_{L'})_-$ is simple for any $1\leq L \leq L'\leq N$. From Proposition \ref{new}, the tensor product $(S_i)_-\otimes (S_j)_-$ is simple for any $i\neq j$. Then all $(S_i)_-$ belong to a category equivalent to $\mathcal{C}_{\ell - 1}$ by Proposition \ref{shift}. Hence the irreducibility of $(S_{L})_-\otimes \cdots \otimes (S_{L'})_-$ follows from the induction hypothesis on $\ell$.

By Corollary \ref{map} we obtain a surjective morphism
$$W
\twoheadrightarrow
\left(S_N\right)_+\otimes\left(S_1\otimes \cdots \otimes S_{N-1}\right)
\otimes \left(S_N\right)_-$$
where $W = \left(S_N\right)_+\otimes\left(\left(S_1\right)_+\otimes \cdots \otimes \left(S_{N-1}\right)_+\right)
\otimes
\left(\left(S_1\right)_-\otimes \cdots \otimes \left(S_{N-1}\right)_-\right)\otimes \left(S_N\right)_-$.
 
We have established above that for every $1\leq i < N$, we have a surjective morphism 
$(S_N)_+\otimes S_i \twoheadrightarrow S_i\otimes (S_N)_+$. Hence we get a sequence of surjective morphisms
\[
\left(S_N\right)_+\otimes\left(S_1\otimes \cdots \otimes S_{N-1}\right)
\twoheadrightarrow
S_1\otimes \left(S_N\right)_+\otimes S_2\otimes\cdots \otimes S_{N-1}
\twoheadrightarrow
\cdots
\twoheadrightarrow
\left(S_1\otimes \cdots \otimes S_{N-1}\right)\otimes (S_N)_+. 
\]
Consequently we get surjective morphisms
$$
\left(S_N\right)_+\otimes\left(\left(S_1\right)_+\otimes \cdots \otimes \left(S_{N-1}\right)_+\right)
\otimes
\left(\left(S_1\right)_-\otimes \cdots \otimes \left(S_{N-1}\right)_-\right)\otimes \left(S_N\right)_-$$
$$\twoheadrightarrow
\left(S_1\otimes \cdots \otimes S_{N-1}\right)\otimes \left(S_N\right)_+
\otimes \left(S_N\right)_-
\twoheadrightarrow
V: = S_1\otimes \cdots \otimes S_N. $$
So $V$ is cyclic since $W$ is cyclic.

Consider the dual module $V^* \cong \left(S_N^*\otimes \cdots \otimes S_1^*\right)$.
By our assumptions, $S_j^*\otimes S_i^* \cong (S_i \otimes S_j)^*$
is simple for every $1\le i<j \le N$. Moreover the modules $S_i^*$ belong to a
category equivalent to $\mathcal{C}_\ell$ by Proposition \ref{shift}. 
So $V^*$ is cyclic.
We can now conclude as in \cite[Section 4.10]{cp}, because a cyclic
module whose dual is cyclic is simple. \qed

\section{Discussions}\label{cinq}

Let us conclude with some comments which are not used in the proof of the main
result of the present paper. 

For simply-laced types, the intermediate Corollary \ref{useqt2} can also be proved by using Nakajima's $q,t$-characters \cite{Nab}. Let us explain this proof since it is related to a nice symmetry property of the corresponding Kazhdan-Lusztig polynomials (a priori, this method can not be extended to the general case since quiver varieties used in \cite{Nab} are not known to exist for the non simply-laced cases).

To start with, let us give some reminders on Nakajima's $q,t$-characters which are certain $t$-deformations
of $q$-characters.

Let $\hat{\Yim}_t = \Yim[Y_{i,a}, V_{i,a}, t^{\pm 1}]_{i\in I,a\in\CC^*}$ which is a $t$-deformation
of $\Yim$. The $V_{i,a}$ are new variables playing the role of the $A_{i,a}^{-1}$ (the $Y_{i,a}$ are denoted by $W_{i,a}$ in \cite{Nab}). 

\noindent A $t$-deformed product $*$ and a bar involution are defined on $\hat{\Yim}_t$ in \cite{Nab}. The bar involution satisfies $\overline{a * b} = \overline{b}*\overline{a}$ for $a,b\in\hat{\Yim}_t$ and $\overline{t} = t^{-1}$. There is a ring morphism $\pi : \hat{\Yim}_t \rightarrow \Yim$ satisfying $\pi(Y_{i,a}) = Y_{i,a}$, $\pi(V_{i,a}) = A_{i,a}^{-1}$, $\pi(t) = 1$ for any $i\in I$, $a\in\CC^*$.

\noindent A monomial $m$ in $\hat{\Yim}_t$ is a product of $Y_{i,a}$, $V_{i,a}$, $t^{\pm 1}$
satisfying $\overline{m} = m$. Let $\mathcal{M}_t$ be the set of these monomials and $B\subset \mathcal{M}_t$ be the set of dominant monomials, that is of $m\in\mathcal{M}_t$ such that $\pi(m)$ is a dominant monomial in $\Yim$.
A dominant monomial $m$ of $\Yim$ is seen as an element of $B$ by the natural identification.

\noindent For $M_1$, $M_2\in\mathcal{M}_t$, we write $M_1\leq M_2$ if $M_1\in M_2 \ZZ[V_{i,a}, t^{\pm 1}]_{i\in I, a\in\CC^*}$. 

\noindent A certain subring $\hat{\mathfrak{K}}_t$ of $(\hat{\Yim}_t,*_t)$ is introduced in \cite{Nab} (it plays a role analog to $\text{Im}(\chi_q)\subset \Yim$). For $i\in I, a\in\CC^*$, there is a unique
$$L_{i,a} \in \hat{\mathfrak{K}}_t \cap \left(Y_{i,a}\left(1 + V_{i,aq} + \sum_{V < V_{i,aq}}\ZZ[t^{\pm 1}] V\right)\right).$$
$L_{i,a}$ is a $t$-analog of $\chi_q(L(Y_{i,a}))$. The existence of $L_{i,a}$ is non trivial and is proved in \cite{Nab} (it can also be proved purely algebraically \cite{her02}). 
We define $t$-analogs of $q$-characters of tensor products of fundamental representations
(or Weyl modules) \cite{Nab}. Let $m\in B$ and set $u_{i,a} = u_{i,a}(\pi(m))$ for $i\in I$, $a\in\CC^*$. Let $E_t(m)\in\hat{\Yim}_t$ equal to
\begin{equation}\label{et}M*\prod_{a\in \left(\CC^*/q^{\ZZ}\right)} \left( \cdots *\left(\prod_{i\in I}\left(L_{i,aq^{-1}}\right)^{\left(* u_{i,aq^{-1}}\right)}\right)*\left(\prod_{i\in I}\left(L_{i,aq^0}\right)^{\left(* u_{i,aq^0}\right)}\right)*\cdots\right),\end{equation}
$$\text{where }M = t^rm\prod_{i\in I,a\in\CC^*}Y_{i,a}^{-u_{i,a}},$$
the $\prod$ denote products for $*$ and $r\in\ZZ$ is set so that $E_t(m)\in m + \sum_{m' < m}\ZZ[t^{\pm 1}]m'$.

For $m\in B$, there exists \cite{Nab} a unique $L_t(m)\in\mathfrak{K}_t$ satisfying
$$\overline{L_t(m)} = L_t(m)\text{ , }L_t(m) \in m + \sum_{m' < m}\ZZ[t^{\pm 1}]m',$$
\begin{equation}\label{kl}E_t(m) = \underset{\left\{m' \in B|m' \leq m\right\}}{\sum}P_{m',m}(t) L_t(m')\text{ with $P_{m',m}(t)\in t^{-1}\ZZ[t^{-1}]\text{ , }P_{m,m}(t) = 1$.}\end{equation}
The $L_t(m)$ are $t$-analogs of $q$-characters of simple modules and the $P_{m',m}(t)$ are analogs of Kazhdan-Lusztig polynomials. The following is an important consequence of the theory of Nakajima's quiver varieties which is proved in \cite{Nab}. We have 
$$\pi\left(L_t(m)\right) = \chi_q\left(L\left(\pi(m)\right)\right)\text{ for $m\in B$.}$$
In particular, $P_{m',m}(1)$ is the multiplicity of $L(\pi(m'))$ in the tensor product of fundamental
representations of highest weight monomial $\pi(m)$.

The $L_t(m)$ can be computed from (\ref{kl}) and from the $E_t(m)$ with an algorithm of Kazhdan-Lusztig type. The $L_{i,a}$, and so the $E_t(m)$, can be computed independently. Hence the $\chi_q(L(m))$ can be computed in principle. In practice, as explained in the introduction of \cite{Nab}, it is difficult to obtain informations on $q$-characters from this complicated algorithm. But theoretical informations can be obtained from it, such as the alternative proof of Corollary \ref{useqt2} for simply-laced types that we sketch now.

We use the notations of Corollary \ref{useqt2}. Let $M' = (M)^{\leq (L - 1)}$ and $M'' = (M)^{\geq L}$. By the defining Formula (\ref{et}), we get $E_t(M) = E_t(M') * E_t(M'')$. By the defining formula of $*$ in \cite{Nab} and Formula (\ref{et}), for $m$ a monomial occurring in $E_t(M'')$, we get $M' * m = m * M' = M' m\in\mathcal{M}_t$. Hence
\begin{equation}\label{etm}E_t(M) = M' E_t\left(M''\right) + \left(E_t\left(M'\right) - M' \right) * E_t\left(M''\right).\end{equation} 
For a monomial $m$ of the form $m = M'' t^r V_{i_1,q^{r_1}}\cdots V_{i_{N(m)},q^{r_{N(m)}}}$, we prove by induction on $N(m)\geq 0$ the following symmetry property :
\begin{equation}\label{sym}[L_t(M)]_{M'm} = [L_t(M'')]_m\text{ and }P_{M' m, M}(t) = P_{m, M''}(t),\end{equation}
where $[\chi]_p\in\ZZ[t^{\pm 1}]$ is the multiplicity of a monomial $p$ in $\chi\in\hat{\Yim}_t$. The property is clear for $N(m) = 0$ since $P_{M,M}(t) = P_{M'',M''}(t) = 1$ and $[L_t(M)]_M = [L_t(M'')]_{M''} = 1$. The inductive step follows from direct computations in the ring $\hat{\Yim}_t$ by using (\ref{kl}) and (\ref{etm}).

\noindent The statement of Corollary \ref{useqt2} is now obtained immediately by applying $\pi$ to (\ref{sym}).

\begin{rem}In the same way we get the following (not used in this paper) 
$$\chi_{q,\leq L}(V)= (M)^{\geq (L + 1)} \chi_{q,\leq L}\left(L\left((M)^{\leq L}\right)\right).$$
Indeed let $M' = (M)^{\leq L}, M'' = (M)^{\geq (L + 1)}$. As above $E_t(M) = E_t(M') * E_t(M'')$ and
$$E_t(M) \in E_t\left(M'\right) M'' + M \sum_{i\in I, r\geq L}\left(V_{i,q^r}\ZZ\left[t^{\pm 1},V_{j,a}\right]_{j\in I, a\in\CC^*}\right).$$
Then, we prove by induction as above and by direct computation that for a monomial $m$ in $M'\ZZ\left[t^{\pm 1},V_{i,q^r}\right]_{i\in I, r < L}$, we have $[L_t(M)]_{M''m} = [L_t(M')]_m$. This means
$$L_t(M) \in M''L_t(M') + M \sum_{i\in I, r\geq L}\left(V_{i,q^r}\ZZ\left[t^{\pm 1},V_{j,a}\right]_{j\in I, a\in\CC^*}\right).$$ 
We are done by applying $\pi$.
\end{rem}

\end{document}